\documentclass[12pt]{article}
\usepackage{amssymb}

\newcommand{\C}{{\mathbb C}}
\newcommand{\LL}{\Lambda}
\newcommand{\EE}{{\mathbb E}}
\newcommand{\PP}{{\rm Prob}}
\newcommand{\Z}{{\mathbb Z}}
\newcommand{\F}{{\mathbb F}}
\newcommand{\A}{{\mathcal A}}
\newcommand{\SUM}{\raisebox{-0.4ex}{\mbox{\Large $\Sigma$}}}
\newcommand{\eps}{\varepsilon}
\newcommand{\B}{{\cal B}}

\newtheorem{theorem}{Theorem}
\newtheorem{proposition}{Proposition}
\newtheorem{lemma}{Lemma}

\title{The structure of critical sets for $\F_p$ arithmetic progressions} 

\author{Ernie Croot}

\begin{document}

\maketitle

\section{Introduction}

Given a function $h : \F_p \times \F_p \times \cdots \times \F_p \to \C$, we define
the usual expectation operator
$$
\EE_{n_1,...,n_k}(h)\ :=\ p^{-k} \SUM_{n_1,...,n_k \in \F_p} h(n_1,...,n_k).
$$
We also define, for $f : \F_p \to \C$, the operator
$$
\LL(f)\ :=\ \EE_{n,d}(f(n)f(n+d)f(n+2d)).
$$
If $f$ were an indicator function for some set $S \subseteq \F_p$, 
this would give a normalized count of the number of three-term progressions in $S$.

In the present paper we establish a new structure theorem for functions 
$f : \F_p \to [0,1]$ that minimize the number of three-term progressions, subject
to a density constraint; and, as a consequence of this result, 
we prove a further structural result, which can also be deduced from 
the work of Green \cite{green}, though only for high densities 
(Green's result only works for densities exceeding $1/\log_*(p)$, though
perhaps his method can be generalized for this particular problem 
to handle lower densities).  

Before stating the theorem, it is worth mentioning that Green and
Sisask \cite{greenolof} have shown that sets of high density (density close
to $1$) that minimize the number of three-term arithmetic progressions,
are the complement of the union of two long arithmetic
progressions (actually, their result is stated in terms of sets that
maximize the number of three-term progressions, but there is a standard
trick to relate this to the minimizing sets).
\bigskip

Our main theorem is now given as follows:

\begin{theorem} \label{main_theorem} Suppose that 
$$
f\ :\ \F_p\ \to\ [0,1]
$$
minimizes $\LL(f)$, subject to the constraint that 
$$
\LL(f)\ \geq\ \theta\ \in\ (0,1].
$$
Then, 
\bigskip

$\bullet$  Let $C(n)$ equal $f(n)$ rounded to the nearest integer, which
is therefore $0$ or $1$.  Then,
$$
\SUM_n |f(n) - C(n)|\ \ll\ p (\log p)^{-2/3}.
$$
So, $f$ must be approximately an indicator function.  Furthermore, we
get the same conclusion if $f$ satisfies $\EE(f) \geq \theta$, and
$\LL(f)$ comes within $O(1/p)$ of the minimal value for this density
constraint.
\bigskip

$\bullet$ There exists a function $r : \F_p \to [0,1]$ 
such that $\EE(r) = \EE(f)$, where $\LL(r)$ is very close to 
the minimal $\LL(f)$, specifically
$$
\LL(r)\ =\ \LL(f)\ +\ O(p^{-1}),
$$
such that if we let, for some $L$,  
$$
S\ :=\ \{ n \in \F_p\ :\ (r*r)(2n) + 2(r*g)(-n) \leq L\},\ 
{\rm where\ } g(n)\ :=\ r(-n/2),
$$ 
then
\begin{equation} \label{rconclusion}
\SUM_n |r(n) - S(n)|\ \ll\ p (\log p)^{-2/3}.
\end{equation}
(Please see subsection \ref{second_part_remark} for an explanation of this
part of the theorem.)

$\bullet$  We have that there exist sets $A$ and $B$ of $\F_p$, with
$|A| > p^{1-o(1)}$ and $|B| > p^{1/2}$, 
such that the set for which $f$ is approximately an indicator function,
is roughly the sumset $A+B$.  More precisely:  If
we let $C(n)$ denote $f$ rounded to the nearest integer, as in 
the first bullet above, then  
$$
\SUM_n |(A*B)(n) - |B| C(n)|\ \ll\ p|B| (\log\log p)^{-2/3}.
$$
Furthermore, we may take $A = C$ and take $B$ to be a certain 
``Bohr neighborhood'' $\B$, which is described in the proof of the theorem.
\end{theorem}

\subsection{A remark on the second part of the theorem} \label{second_part_remark}
 
By (\ref{rconclusion}) we see that $r$ is nearly an indicator function for
the set $S$.  Let us suppose, for the purposes of discussion, 
that it is exactly an indicator function for some set, and let $R$ be this set.  
Note that $R$ and $S$ must have small symmetic difference.

When does an $n \in \F_p$ belong to the set $S$?  To decide this, given $n$,
we let $N_1$ be the number of pairs $(x,y) \in R \times R$ such that 
$n,x,y$ forms an arithmetic progression; we let $N_2$ be the number of 
pairs $(x,y) \in R \times R$ such that $x,n,y$ is an arithmetic progression; and,
we let $N_3$ be the number of such ordered pairs where $x,y,n$ is an 
arithmetic progression.  For $n$ to belong to $S$, we must have that
$$
N_1 + N_2 + N_3\ \leq\ L.
$$ 

Since $S$ and $R$ have small symmetric difference, we see that 
our $r$ can be thought of as enjoying a ``local minimal'' property:  
Not only does $r$ minimize $\LL(r)$ up to an error $O(1/p)$, 
subject to $\EE(r) \geq \theta$, but we can easily decide whether $n \in \F_p$
belongs to $R$ or not, simply by checking to see how 
many progressions pass through the point $n$, with the 
other two end-points in $R$.  If this count is
small enough, then $n$ likely belongs to $R$ (though it certainly
belongs to $S$); but, if the count is
large, $n$ likely does not belong to $R$.

The most difficult part of this proof that $R$ and $S$ are nearly
the same, is handling those $n$ where $N_1 + N_2 + N_3$ exactly equals $L$.

\subsection{Remarks on the third part of the theorem}

One reason to believe the third bullet above
is that from the second bullet we expect that $f$ is an indicator
function for a level set of a ``smooth function'' $(r*r)(2n) + 2(r*g)(-n)$;
and, as is well known, such level sets must be approximately the union
of a bunch of translates of a Bohr neighborhood of the function, at least
when their density is large enough.  
\bigskip

It should be remarked that sumsets are quite special structures, as
are smooth functions of the type $(r*r)(2n) + 2(r*g)(-n)$, and only
a vanishingly small proportion of the subsets of $\F_p$ are sumsets
or form the support of a smooth function; so,
the third bullet is saying something fairly non-trivial about our 
minimal $f$.  

Also, there are loads of other consequences that one 
can deduce from the third
bullet.  One of these is that, upon decomposing the Bohr neighborhood $\B$
into a union of arithmetic progressions, 
one can deduce that $C$ is essentially the union of a ``small number'' of  
somewhat ``long'' arithmetic progressions (``small number'' can mean
a power of $p$, say $p^c$, where $c < 1$), all having the same common difference.

\section{Proof of Theorem \ref{main_theorem}}

The proof of this structure theorem depends on a certain function 
$r_3$, which we presently define.
\bigskip

\noindent {\bf Definition.}  Given a subset $S$ of a group $G$,
we let $r_3(S)$ denote the size of the largest subset of $S$ 
free of solutions to $x+y = 2z$, $x \neq y$.  In all the uses of $r_3$ in
the present paper, $G = \Z$ and $S = [N] := \{1,2,...,N\}$, for various
different values of $N$.
\bigskip

\noindent Bourgain \cite{bourgain} has recently shown that
\begin{equation} \label{bourgainsresult}
r_3([N])\ \ll\ N (\log N)^{-2/3},
\end{equation}
and from a result of Behrend \cite{behrend}, we know that for $N$ sufficiently large,
$$
r_3([N])\ >\ N \exp( - c \sqrt{\log N}),
$$
for a certain constant $c > 0$.

\subsection{Proof of the first part of Theorem \ref{main_theorem}} 
\label{first_part}

For this part we will begin by assuming that 
$\EE(f) > \kappa p (\log p)^{-2/3}$, for as large a $\kappa > 0$ as
we might happen to need, since this part of the theorem is 
trivially true otherwise.

Here we will first show that the minimal $f$ is well-approximated by
an indicator function; actually, we will prove even more -- we will
show that if $\LL(f)$ comes within $O(p^{-1})$ of this smallest value,
subject to the density constraint $\EE(f) > \theta$,
then $f$ must be approximately an indicator function.
To do this, we will require the 
following proposition, proved in subsection \ref{prop_subsection}.

\begin{proposition}  \label{level_prop} 
Suppose that $A$ and $B$ are disjoint subsets of $\F_p$, such that 
$f : \F_p \to [0,1]$ has the property  
$$
{\rm for\ } n\in A,\ f(n)\ \leq\ 1- \eps,\ 0 < \eps\ <\ 1/3,
$$
and suppose that
$$
{\rm support}(f)\ =\ A \cup B.
$$  
Then, for $\beta > 0$ satisfying 
$$
\eps \beta\ \geq\ p^{-1/2} \log p,
$$
there exists a function $g : \F_p \to [0,1]$ such that 
$$
\EE(g)\ \geq\ \EE(f),
$$
and yet
$$
\LL(g)\ <\ \LL(f) + 2\beta - \eps^2 p^{-2} W_0/4 + O(p^{-1}),
$$
where
$$
W_0\ :=\ \SUM_{a,a+d,a+2d \in A} f(a)f(a+d)f(a+2d).
$$
\end{proposition}

We also will require the following quantitative version of Varnavides's
theorem \cite{varnavides}.

\begin{lemma} \label{quantitative_varnavides}  
If $S \subseteq \F_p$ satisfies $|S| \geq 2 (r_3(N)/N)p$, we will have
for any $2 \leq N \leq p$ that 
$$
\LL(S)\ \geq\ {2r_3([N]) \over N^3 + O(N^2)}.
$$
\end{lemma}

\noindent {\bf Proof of the Lemma.}  The proof of this lemma is via
some easy averaging:  We let $\A_N$ denote the set of all arithmetic
progressions $A \subseteq \F_p$ having length $N$.  These arithmetic
progressions are to be identified by ordered pairs $(a,d)$, $d \neq 0$,
where $a$ is the first term in the progression, and where $d$ is the common 
difference.  Note that this means we ``double count'' arithmetic progressions
in that the progression $a,a+d,a+2d,...,a+kd$ is distinct from 
$a+kd, a+(k-1)d,...,a$.  

It is easy to check that each sequence $a,a+d,a+2d$, $d \neq 0$ 
is contained in exactly $N^2/2 + O(N)$ of these $A \in \A_N$:  
We have that each three-term
progression is contained in the same number of $A \in \A_N$, and each 
$A \in \A_N$ contains $N^2/2 + O(N)$ three-term progressions; hence,
if $P$ denotes the number of $A \in \A_N$ containing a particular 
sequence $a,a+d,a+2d$, we have since there are $p(p-1)$ non-trivial 
progressions in $\F_p$, that
$$
p(p-1) P\ =\ |\A_N| (N^2/2 + O(N)),
$$
whence $P = N^2/2 + O(N)$.

So, if we let $T_3(X)$ denote the number of sequences
$a,a+d,a+2d \in X$, $d \neq 0$, we have that 
\begin{equation} \label{T3}
T_3(S)\ =\ \left ( N^2/2 + O(N) \right )^{-1} \SUM_{A \in \A_N} T_3(A \cap S).
\end{equation}
Next, we need a lower bound on how many $A \in \A_N$ satisfy 
$|A \cap S| \geq r_3(N)$:  First, note that for each $d \in \F_p$,
$d \neq 0$, there are exactly $N$ arithmetic progressions $A \in \A_N$ having 
common difference $d$ that contain a particular point $a \in \F_p$.  So, 
$$
\SUM_{A \in \A_N} |A \cap S|\ =\ \SUM_{s \in S} \SUM_{d \in \F_p \atop d \neq 0} N
\ =\ (p-1) N |S|.
$$
Let $Y$ be the number of $A \in A_N$ for which $|A \cap S| > r_3(N)$.
Then, we have
$$
(|\A_N| - Y) r_3(N) + Y N\ \geq\ (p-1)N |S|,
$$
which implies 
$$
Y\ \geq\ {(p-1) N |S| - |\A_N| r_3(N) \over N - r_3(N)}\ \geq\ (p-1) |S| - |\A_N| (r_3(N)/N).
$$

For each of these $Y$ progressions $A \in \A_N$ we will have that 
$T_3(A \cap S) \geq 1$; and so, we deduce from (\ref{T3}) that 
$$
T_3(S)\ \geq\ {(p-1)|S| - |\A_N| (r_3(N)/N) \over N^2/2 + O(N) }.
$$
Using the easy to see fact that $|\A_N| = p(p-1)$, we deduce that if 
$$
|S|\ >\ 2 (r_3(N)/N) p,
$$
then
$$
T_3(S)\ \geq\ {2 p^2 (r_3(N)/N) \over N^2 + O(N)}.
$$
The lemma easily follows on rephrasing this in terms of $\LL(S)$.
\hfill $\blacksquare$
\bigskip

Now we let
$$
A\ :=\ \{ n \in \F_p\ :\ f(n) \in [\eps, 1- \eps]\},
$$
where $\eps > 0$ will be determined later.  In order for
$f$ to be minimal, from Proposition \ref{level_prop} we 
deduce that we must have that if  
$\eps \beta = p^{-1/2} \log p$, then  
$$
\beta\ \geq\ \eps^2 p^{-2} W_0/8 + O(1/p).
$$
So, since we trivially have that
$$
W_0\ \geq\ \eps^3 p^2 \Lambda(A),
$$
it follows that  
\begin{equation} \label{LLA}
\LL(A)\ \leq\ 8\eps^{-6} p^{-1/2} \log p.
\end{equation}

We would like to now apply Lemma \ref{quantitative_varnavides} to this,
but in order to do so, we must solve for $N$ such that 
$$
|A|\ >\ 2r_3(N) p/N.
$$
To this end, we require the bound (\ref{bourgainsresult}) of Bourgain,
which implies that if we let 
$$
N\ =\ \exp ( c(p/|A|)^{3/2} )\ <\ p,\ {\rm since\ } 
|A| > \kappa p (\log p)^{-2/3},
$$
then we will have that 
$$
|A|\ >\ p (\log N)^{-2/3}\ >\ 2r_3(N)p/N,
$$
as we require.

From this it follows from Lemma \ref{quantitative_varnavides} that
$$
\LL(A)\ >\ r_3(N)/N^3\ >\ 1/N^3\ >\ \exp(-3c (p/|A|)^{3/2}).  
$$
It follows now from (\ref{LLA}) that
$$
|A|\ \ll\ p \log^{-2/3} (\eps^{12} p),\ {\rm for\ } \eps > p^{-1/12} \log p.
$$
So, if we let $C$ be the function $f$ rounded to the nearest integer
(which will be either $0$ or $1$), then for $n \in A$ we will have 
$|f(n) - C(n)| \leq 1$, while for all other $n$ we will have 
$|f(n) - C(n)| \leq \eps$.  It follows that 
$$
\SUM_n |f(n) - C(n)|\ \ll\ (\eps + (\log \eps^{12} p)^{-2/3})p,
\ {\rm for\ } \eps\ >\ p^{-1/12}\log p.
$$
Choosing $\eps = (\log p)^{-2/3}$, we deduce that this sum is
$O(p (\log p)^{-2/3})$, just as in Bourgain's theorem 
(\ref{bourgainsresult}).  This completes the proof of the first part 
of our theorem.

\subsection{Proof of the second part of Theorem \ref{main_theorem}}

Given a function $h : \F_p \to [0,1]$, we let 
$$
h_2(n)\ :=\ h(-n/2),
$$
and then we define
$$
F_h(n)\ :=\ (h*h)(2n) + (h*h_2)(-n).
$$

In order to proceed further, we will require the following proposition.

\begin{proposition} \label{forcing_prop}
Fix $A \subseteq \F_p$, and associate to each $a \in A$ a real number 
$w_a \in [0,1]$.  Among all functions $h : \F_p \to [0,1]$ satisfying 
$$
\EE(h)\ =\ \gamma\ >\ \SUM_{a \in A} w_a,
$$
those which minimize $\LL(h)$ have the property that there exists 
$L > 0$ such that  
$$
{\rm for\ } n \in \F_p \setminus A,\ 
h(n)\ =\ \left \{ \begin{array}{rl} 1,\ & {\rm if\ } F_h(n) < L; \\
0,\ & {\rm if\ } F_h(n) > L. \end{array} \right .
$$
\end{proposition}
The proof of this proposition can be found in subsection
\ref{forcing_prop_section}.  

We will use this proposition to construct a sequence of sets
$$
A_1,\ A_2, ...\ \subseteq\ \F_p,
$$
and a sequence of functions
$$
r_1,\ r_2,\ ...,\ :\ \F_p\ \to\ [0,1],\ {\rm and\ all\ }\EE(r_i)\ =\ \EE(f).
$$
such that the following all hold.
\bigskip

$\bullet$ First, $|A_1| = 2$ and $A_{i+1} = A_i \cup \{x_{i+1}, y_{i+1}\}$; 

$\bullet$ second, $\LL(r_i) \leq \LL(f) + 5 p^{-2}|A_i|$; 

$\bullet$ third, given particular fixed values for $r_i(n)$ on 
$A_i$, we have that $r_i$ minimizes $\LL(r_i)$, subject to the density
constraint $\EE(r_i) = \EE(f)$;

$\bullet$ and finally, for each $n \in A_i$, $r_i(n) \in [1/4, 3/4]$.
\bigskip

\noindent Clearly, this process cannot continue past 
the $\lfloor p/2\rfloor$th iteration, as the sets $A_i$ grow
by two elements after each iteration.  Furthermore, we will show
that whenever the process {\it does} terminate (which it will in
either case 1 or case 2 below), we will be left with
a function $r : \F_p \to [0,1]$ satisfying the conclusion 
in the second bullet of Theorem \ref{main_theorem}.

For the time being, let us suppose that these sequences can 
be constructed as claimed:  Suppose we have
constructed $A_i$; we will now show how to construct $A_{i+1}$.
To this end, we apply Proposition \ref{forcing_prop} with
$A = A_i$ (in the case $i=0$ we let $A$ be the empty set), and then
we deduce that for some $L > 0$, $r_i(n) = 1$ (we use $r_0(n) := f(n)$) 
for $F_{r_i}(n) < L$ and $r_i(n) = 0$ for $F_{r_i}(n) > L$.  
We furthermore apply the already-proved first part of 
Theorem \ref{main_theorem} from subsection \ref{first_part}, and deduce
that since $\LL(r_i) \leq \LL(f) + 10i p^{-2}$,
$$
\SUM_n |r_i(n) - C(n)|\ \ll\ p (\log p)^{-2/3},
$$
where $C(n)$ is $r_i(n)$ rounded to the nearest integer.
\bigskip

Now, if there are two distinct places $x,y \in \F_p \setminus A_i$ for which
$$
r_i(x), r_i(y)\ \in\ [1/4,3/4],
$$
then we just let 
$$
A_{i+1}\ :=\ \{x,y\},\ {\rm and\ }r_{i+1}\ :=\ r_i.
$$  

So suppose that there are no such $x$ and $y$; there are three possibilities
to consider.

\subsubsection{Case 1:  $r_i(n) \geq 1/2$ for all $n \in \F_p \setminus A_i$ 
where $F_{r_i}(n) = L$.}  

Note that we include in this case the possibility that there are no 
$n \in \F_p \setminus A_i$ such that $F_{r_i}(n) = L$.

If we are in this case, then it means that 
$$
\SUM_{n\in \F_p \setminus A_i\ :\ F_{r_i}(n) = L} |r_i(n) - 1|\ \ll\ 
p (\log p)^{-2/3};
$$
and so, $r_i(n)$ is very close to $1$ at most places $n \in \F_p \setminus A_i$ 
where $F_{r_i}(n) = L$.  It follows that if we were to let $S$ be the 
set of all $n \in \F_p \setminus A_i$ with $F_{r_i}(n) \leq L$, then 
$$
\SUM_{n \in S} |r_i(n) - 1|\ \ll\ p (\log p)^{-2/3}.
$$
In order to extend this sum to all $n \in \F_p$, we will need to 
show that $|A_i|$ cannot be too big.  Basically, we will show that if
it is, then $\LL(f)$ could not be minimal.  

To see this last point, we apply Proposition \ref{level_prop}, using 
$A := A_i$, $h := r_i$, $\eps = 1/4$, and $\eps \beta = p^{-1/2} \log p$,
and we deduce that there exists 
$$
g\ :\ \F_p\ \to\ [0,1],\ \EE(g)\ \geq\ \EE(r_i),
$$
and yet 
$$
\LL(g)\ \leq\ \LL(f) + 8 p^{-1/2}\log p - p^{-2} W_0/64 + O(p^{-1}),
$$
where
$$
W_0\ :=\ \SUM_{a,a+d,a+2d \in A} r_i(a)r_i(a+d) r_i(a+2d)\ \geq\ 
4^{-3} p^2 \LL(A).
$$ 
We wish to apply Lemma \ref{quantitative_varnavides}:  First, let 
$$
N\ =\ \exp ( c(p/|A_i|)^{3/2} )\ <\ p,
$$
such that from (\ref{bourgainsresult}) we deduce that
$$
|A_i|\ >\ 2p (\log N)^{-2/3}\ >\ 2r_3([N])p/N,
$$
as we require.

From this it follows now from Lemma \ref{quantitative_varnavides} that
\begin{equation} \label{LA}
\LL(A_i)\ \geq\ {2r_3(N) \over N^3 + O(N^2)}\ >\ 1/N^3\ \gg\ 
\exp(-3c(p/|A_i|)^{3/2}),
\end{equation}
for $N$ sufficiently large.

In order for $f$ to minimize $\LL(f)$, we must have that 
$$
8 p^{-1/2} \log p\ =\ 2\beta\ \geq\ p^{-2} W_0/64 + O(1/p);
$$
so, ignoring the $O(1/p)$, we see that
$$
\LL(A_i)\ \leq\ 4^3 p^{-2} W_0\ \leq\ 4^8 p^{-1/2} \log p. 
$$
It follows from this and (\ref{LA}) that
$$
|A_i|\ \ll\ p (\log p)^{-2/3}, 
$$
as claimed.  It follows that if we extend $S$ to be the 
set of all $n$ where $F_{r_i}(n) \leq L$, then  
$$
\SUM_n |r_i(n) - S(n)|\ \ll\ p (\log p)^{-2/3},
$$
and the second bullet of Theorem \ref{main_theorem} is proved upon
setting $r = r_i$.

\subsubsection{Case 2: $r_i(n) < 1/2$ for all $n \in \F_p \setminus A_i$ 
where $F_{r_i}(n) = L$.}

If we are in this case, then it means that 
$$
\SUM_{n\in \F_p \setminus A_i\ :\ F_{r_i}(n) = L} r_i(n)\ \ll\ 
p (\log p)^{-2/3};
$$
and so, if we let $L' = L - \delta$, for small enough $\delta > 0$,
then we will have that for $n \in \F_p \setminus A_i$, 
$r_i(n) = 1$ for $F_{r_i}(n) \leq L'$, while 
$r_i(n)$ is usually near $0$ when $F_{r_i}(n) > L'$.  
It follows then that 
if we let $S$ be the set of $n \in \F_p\setminus A_i$ where 
$F_{r_i}(n) \leq L'$, then 
$$
\SUM_{n \in S} |r_i(n) - 1|\ \ll\ p (\log p)^{-2/3},
$$

We wish to extend this to where $S$ is the set of all $n$ satisfying 
$F_{r_i}(n) \leq L'$, by showing that $|A_i|$ cannot be too big, 
and we proceed exactly the same
way as in Case 1 above.  We then deduce that, upon redefining $S$ in 
this way, that
$$
\SUM_n |r_i(n) - S(n)|\ \ll\ p (\log p)^{-2/3},
$$
and again this proves the second bullet of Theorem \ref{main_theorem}
upon setting $r = r_i$.

\subsubsection{Case 3:  There exists $x,y \in \F_p \setminus A_i$ where 
$F_{r_i}(x) = F_{r_i}(y) = L$, and $r_i(x) < 1/2 < r_i(y)$.}

In order to decide what to do in this case, we will require the following
basic fact, which is an immediate consequence of the 
formula for $\LL(h_3)$ in the proof of Proposition \ref{forcing_prop}
in section \ref{forcing_prop_section} in equation
(\ref{Lf2}):  We have that if 
we let 
$$
r_{i+1}(n)\ :=\  \left \{ \begin{array}{rl} r_i(n),\ & {\rm if\ }n \neq x,y; \\
(r_i(x) + r_i(y))/2,\ & {\rm if\ } n = x\ {\rm or\ }y, \end{array}\right.
$$
then  
\begin{eqnarray*}
\LL(r_{i+1})\ &\leq&\ \LL(r_i)\ +\ p^{-2} (r_{i+1}(x) - r_i(x))F_{r_i}(x) + 
p^{-2} (r_{i+1}(y) - r_i(y)) F_{r_i}(y) + 10 p^{-2} \\
&\leq&\ \LL(r_i) + 10 p^{-2}.
\end{eqnarray*}
So, when we are in this case, we just let 
$$
A_{i+1}\ :=\ \{x,y\},
$$
and so the properties of $A_{i+1}, r_{i+1}$ that we require all  
hold.

\subsection{Proof of the third part of Theorem \ref{main_theorem}}

We assume for this part of the proof of our theorem that 
$\theta > (\log\log p)^{-2/3}$, since our problem is trivial otherwise.

We now prove the third bullet of Theorem \ref{main_theorem}. 
To this end, we let 
$$
f_3(n)\ :=\ (f*\mu)(n),
$$
where $\mu$ is defined as follows:  First, we locate the places 
$b_1,...,b_t$ where the Fourier transform
$$
|\hat f(b_i)|\ >\ \eps_0 p,
$$
where $\eps_0 > 0$ will be decided later, and then we define the Bohr 
neighborhood $\B$ to be all those $n \in \F_p$ where
$$
|| b_i n/p||\ <\ \eps_0,\ {\rm for\ all\ }i=1,...,t.
$$
Finally, we just let $\mu(n) = 1/|\B|$ if $n \in \B$, and $\mu(n) = 0$
otherwise.

Our goal now will be to show that 
\begin{equation} \label{ourgoal}
\SUM_n |f_3(n) - f(n)|\ \ll\  p(\log\log p)^{-2/3},
\end{equation}
for this will imply the third bullet of Theorem \ref{main_theorem}
holds:  To see this, note that from the already-proved first bullet,
we know that if we let $C(n)$ be $f(n)$ rounded to the nearest integer,
then 
\begin{eqnarray*}
\SUM_n | |\B|^{-1} (C*\B)(n) - C(n)|\ &=&\ 
\SUM_n | |\B|^{-1} (f*\B)(n) - f(n)| + O(p (\log p)^{-2/3}) \\
&=&\ \SUM_n |f_3(n) - f(n)| + O(p (\log p)^{-2/3})\\
&\ll&\ p(\log\log p)^{-2/3},
\end{eqnarray*}
which is just what the third bullet claims.
\bigskip

Now we show that (\ref{ourgoal}) holds:  First note that Parseval gives
$$
t\ \leq\ \theta \eps_0^{-2};
$$
and the following standard lemma tells us that our Bohr neighborhood
is ``large''.

\begin{lemma}  We have that 
$$
|\B|\ \geq\ (\eps_0 + O(1/p))^t p. 
$$
\end{lemma} 

\noindent {\bf Proof of the lemma.}  For $i=1,2,...,t$, we let 
$$
\alpha_i(x)\ :=\ (\eps_0 p + 1)^{-1} 
\left (\SUM_{||b_i n/p|| < \eps_0/2} e^{2\pi i n x/p} \right )^2
$$
We note that $\alpha_i(x)$ is always a non-negative real for all real
numbers $x$, and $\alpha_i$ is the Fourier transform of a function
$\beta_i : \F_p \to [0,1]$.  Furthermore, 
$$
|\alpha_i(0)|\ =\ \eps_0 p + O(1).
$$
Now letting 
$$
\beta(n)\ :=\ (\beta_1 \cdots \beta_t)(n),
$$
we find that $\beta : \F_p \to [0,1]$, and has support contained within $\B$.
So,
\begin{eqnarray*}
|\B|\ \geq\ \hat \beta(0)\ &=&\ p^{-t+1} (\hat \beta_1 * \hat \beta_2 * 
\cdots * \hat \beta_t)(0)\\
&=&\ p^{-t+1} (\alpha_1 * \alpha_2 * \cdots * \alpha_t)(0) \\
&\geq&\ p^{-t+1} \alpha_1(0)\cdots \alpha_t(0) \\
&\geq&\ (\eps_0 + O(1/p))^t p.
\end{eqnarray*}

\hfill $\blacksquare$
\bigskip

Now, from the easy-to-check fact that
$$
||\hat f_3(a) - \hat f(a)||_\infty\ =\ ||\hat f(a)(1 - \hat \mu(a))||_\infty 
\leq\ \eps_0 p,
$$
we easily deduce, via standard arguments (Parseval and Cauchy-Schwarz) 
that
\begin{eqnarray*}
\LL(f_3)\ =\ p^{-3} \SUM_a \hat f_3(a)^2 \hat f_3(-2a) 
&=&\ p^{-3} \SUM_a \hat f(a)^2 \hat f(-2a)\ +\ E \\
&=&\ \LL(f)\ +\ E,
\end{eqnarray*}
where the ``error'' $E$ satisfies
$$
|E|\ \leq\ 10 \eps_0.
$$
\bigskip

Now let $A$ be all those $n \in \F_p$ for which 
$$
f_3(n)\ \in\ [\eps_1, 1-\eps_1].
$$
Then, we have that 
$$
W_0\ :=\ \SUM_{a,a+d,a+2d \in A} f_3(a)f_3(a+d)f_3(a+2d)\ \geq\ 
\eps_1^3 p^2 \LL(A).
$$
In order to apply Lemma \ref{quantitative_varnavides} to this, we
let 
$$
N\ =\ \exp(c (p/|A|)^{3/2})\ <\ p,
$$
so that from (\ref{bourgainsresult}) we deduce that
$$
|A|\ >\ p (\log N)^{-2/3}\ >\ 2r_3(N)p/N,
$$
as we require.

From this it follows now from Lemma \ref{quantitative_varnavides} 
that 
$$
\LL(A)\ \geq\ {2r_3(N) \over N^3 + O(N^2)}\ >\ 1/N^3\ \gg\ 
\exp(-3(2p/|A|)^{3/2}),
$$ 
for $N$ sufficiently large.

In order for $\LL(f)$ to be minimal, we must have that 
$$
\LL(f)\ \leq\ \LL(f_3)\ \leq\ \LL(f) + 2\beta + 10 \eps_0 - 
\eps_1^2 p^{-2}W_0/4 + O(1/p).
$$
Setting $\beta = 5\eps_0$ we must have
$$
20\eps_0\ \geq\ \eps_1^2 p^{-2}W_0/2 + O(1/p)\ \geq\ \eps_1^5 \LL(A)/2
+ O(1/p);
$$
and so,
$$
\LL(A)\ \leq\ 80 \eps_0 \eps_1^{-5} + O(1/p).
$$
Combining this with our lower bound for $\LL(A)$ above, we deduce that
$$
|A|\ \ll\ p (\log \eps_1^5 \eps_0^{-1})^{-2/3}.
$$

It now follows that if $C(n)$ is $f_3(n)$ rounded to the nearest integer,
then
\begin{eqnarray*}
\SUM_n |f_3(n) - C(n)|\ &\leq&\ \SUM_{n \in A} 1/2\ +\ \SUM_{n \in \F_p 
\setminus A} \eps_1 \\
&\ll&\ p (\log \eps_1^5 \eps_0^{-1})^{-2/3} + \eps_1 p.
\end{eqnarray*}

Now we will set
$$
\eps_0\ :=\ \sqrt{\theta \log\log p/\log p},\ {\rm and\ } 
\eps_1\ :=\ (\log\log p)^{-2/3},
$$
which will give
$$
|\B|\ >\ p^{1/2},
$$
and then our sum on $|f_3(n) - C(n)|$ will be at most 
$$
\SUM_n |f_3(n) - C(n)|\ \ll\ p (\log\log p)^{-2/3},
$$
which completes the proof of Theorem \ref{main_theorem}.

\subsection{Proof of Proposition \ref{level_prop}}  \label{prop_subsection}

\subsubsection{Technical lemmas needed for the proof of the Proposition}

We will need to assemble some lemmas to prove this proposition.  We begin
with the following standard fact:

\begin{lemma} \label{complement_lemma}  Suppose that $S \subseteq \F_p$ 
satisfies $|S| = \alpha p$.  Let $T$ denote the complement of $S$.  Then, we have
that 
$$
\LL(S) + \LL(T)\ =\ 1 - 3\alpha + 3\alpha^2.
$$
\end{lemma}

\noindent {\bf Proof of the lemma.}  One way to prove this is via Fourier analysis:  We
have that 
$$
\LL(S) + \LL(T)\ =\ p^{-3} \SUM_a (\hat S(a)^2 \hat S(-2a) + \hat T(a)^2 \hat T(-2a)).
$$
Since $\hat S(a) = -\hat T(a)$ for $a \neq 0$, we have that all the terms except for
$a = 0$ vanish.  So,
$$
\LL(S) + \LL(T)\ =\ p^{-3} (\hat S(0)^3 + \hat T(0)^3)\ =\ \alpha^3 + (1-\alpha^3)\ =\ 
1 - 3 \alpha + 3 \alpha^2. 
$$
\hfill $\blacksquare$
\bigskip

From this lemma, one can deduce the following corollary, which we state as another
lemma:  

\begin{lemma} \label{few_ap_lemma}  
For $\alpha > 2/3$ we have that there exists a set $S \subseteq \F_p$ 
satisfying $|S| = \lfloor \alpha p \rfloor$, and 
$$
\LL(S)\ \leq\ \alpha^3 (1 - (1-\alpha)^2/2) + O(1/p).
$$
\end{lemma}
\bigskip

\noindent {\bf Proof of the Lemma.} Let $\beta = 1 - \alpha < 1/3$, 
and then let $S$ just be the arithmetic progression
$\{0,1,...,\lfloor \alpha p \rfloor-1 \}$, and then let $T$ be the
complement of $S$, which is also just an arithmetic progression.  
It is easy to check that 
$$
\LL(T)\ =\ |T|^2/2p^2 + O(|T|/p^2)\ =\ \beta^2/2 + O(1/p),
$$
as the solutions to $x+y=2z$, $x,y,z \in T$ are exactly those ordered pairs 
$(x,z) \in T \times T$ of the same parity.  

Applying Lemma \ref{complement_lemma} to this set $T$, we find that 
\begin{eqnarray*}
\LL(S)\ &=&\ (1 - 3\beta + 3\beta^2) - \beta^2/2 + O(1/p)\\
&=&\ 1 - 3 \beta + 5\beta^2/2 + O(1/p) \\
&<&\ (1 - \beta)^3 (1 - \beta^2/2) + O(1/p),
\end{eqnarray*}
as claimed.
\hfill $\blacksquare$

\subsubsection{Body of the proof of Proposition \ref{level_prop}}

We will define the function $g : \F_p \to [0,1]$ such that 
$$
{\rm support}(g)\ \subseteq\ A \cup B,
$$
where
$$
{\rm for\ }n \in B,\ g(n)\ =\ f(n),
$$
but on the set $A$, the funciton $g$ will be different from $f$:   
Basically, we let $S$ be the set produced by Lemma \ref{few_ap_lemma} with 
$\alpha = 1 - \eps$, then take $T$ to be a random
translate and dilate of $S$, say 
$$
T\ :=\ m.S + t\ =\ \{ms + t\ :\ s \in S\}.
$$
Then, we let 
$$
{\rm for\ }n \in A,\ g(n)\ =\ (1-\eps)^{-1} f(n) T(n).
$$
Note that this is $\leq 1$, because we know $f(n) \leq 1 - \eps$ on $A$.

We will show that, so long as there are ``enough'' three-term progressions
lying in $A$, this new function $g$ will have the property that $\LL(g)$ 
is much smaller than $\LL(f)$.  To this end, we consider three types of arithmetic
progressions that give rise to the counts $\LL(f)$ and $\LL(g)$:  Those progressions
that pass through both $A$ and $B$ (say one point in $A$ and two in $B$; or 
two in $A$ and one in $B$); those that lie entirely within $A$; and those that
lie entirely within $B$.  

The contribution to $\LL(g)$ of those arithmetic progressions lying entirely within 
$B$ is the same as the contribution to $\LL(f)$.  So, we don't need to account for
these when trying to prove our upper bound on $\LL(g)$; and therefore there are only two
non-trivial cases that we need to work out:
\bigskip

\noindent {\bf Case 1 (all three points in $A$).}
\bigskip

Define the random variable
$$
Z_0\ :=\ \SUM_{a,a+d,a+2d \in A} g(a)g(a+d)g(a+2d),
$$
and let $W_0$ be the analogous sum but with $g$ replaced by $f$.  
We note that if we only consider those terms with $d\neq 0$, we lose
at most $O(p)$ in estimating $Z_0$.  

We have that
\begin{eqnarray*}
\EE(Z_0)\ &=&\ \SUM_{a,a+d,a+2d \in A \atop d \neq 0} 
\EE(g(a) g(a+d) g(a+2d)) + O(p) \\ 
&=&\ p^{-2} (1-\eps)^{-3} \SUM_{a,a+d,a+2d \in A \atop d \neq 0} 
f(a)f(a+d)f(a+2d) 
\SUM_{m,t \in \F_p \atop a,a+d,a+2d \in m.S + t} 1 + O(p) \\
&=&\ p^{-2} (1-\eps)^{-3} \SUM_{a,a+d,a+2d \in A \atop d \neq 0} 
\SUM_{b,b+d',b+2d' \in S} \\
&&\hskip1in \SUM_{m,t \in \F_p \atop mb + t = a,\ m(b+d') + t = a+d} 
f(a)f(a+d)f(a+2d) + O(p).
\end{eqnarray*}
To estimate this inner sum, we note that the contribution of those
terms with $d' = 0$ is $0$; and, when $d' \neq 0$, 
we get a contribution of $f(a)f(a+d)f(a+2d)$ to just
the inner sum, because there is only one pair $m,t$ which works.  
Thus, we deduce from this and Lemma \ref{few_ap_lemma} that
\begin{eqnarray*}
\EE(Z_0)\ &=&\ p^{-2}(1-\eps)^{-3} 
\SUM_{b,b+d',b+2d' \in S \atop a,a+d, a+2d \in A} 
f(a)f(a+d)f(a+2d)\ +\ O(p) \\
&=&\ (1-\eps)^{-3} \LL(S) W_0\ +\ O(p) \\
&<&\ (1 - \eps^2/2) W_0\ +\ O(p).
\end{eqnarray*}
\bigskip

\noindent {\bf Case 2 (at least one point in $A$, and at least one in $B$).}
\bigskip

Define the random variables
\begin{eqnarray*}
Z_1\ &:=&\ \SUM_{a,a+d \in A \atop a+2d \in B} g(a)g(a+d)g(a+2d) \\
Z_2\ &:=&\ \SUM_{a,a+2d \in A \atop a + d \in B} g(a) g(a+d)g(a+2d) \\
Z_3\ &:=&\ \SUM_{a+d,a+2d \in A \atop a \in B} g(a) g(a+d)g(a+2d) \\
Z_4\ &:=&\ \SUM_{a \in A \atop a+d,a+2d \in B} g(a) g(a+d)g(a+2d) \\
Z_5\ &:=&\ \SUM_{a+d \in A \atop a,a+2d \in B} g(a) g(a+d)g(a+2d) \\
Z_6\ &:=&\ \SUM_{a+2d \in A \atop a,a+d \in B} g(a) g(a+d)g(a+2d).
\end{eqnarray*}
Also, let $W_1,...,W_6$ be the analogous constants with $g$ replaced by $f$
(note that these are not random variables).

We will now compute the expectations of these random variables; though, we will
not do all of these here, and instead will just 
work it out for $Z_1$, as showing it for all the others can be done in 
exactly the same way, and leads to the same bounds.

We have that
$$
\EE(Z_1)\ =\ \SUM_{a+2d \in B} f(a+2d) \SUM_{a,a+d \in A} \EE(g(a)g(a+d)).
$$
To evaluate this last expectation, let us suppose that $a + 2d\in B$ and
$a,a+d \in A$, where $d \neq 0$ (if $d=0$ then we would have that $a$ 
lies both in $A$ and $B$, which is impossible).  Then, given any pair of 
distinct elements $x,y \in S$, there exists a unique pair $(m,t) \in \F_p \times \F_p$
such that 
$$
mx + t\ =\ a\ \ {\rm and\ \ } my + t\ =\ b.
$$
So, the probability that 
$$
g(a)g(a+d)\ =\ (1-\eps)^{-2} f(a)f(a+d),
$$
given $a+2d \in B$, $a,a+d \in A$, is $1/p^2$ times the number of 
ordered pairs $(x,y)$ of distinct elements of $S$, which is 
$|S|(|S|-1)$.  Note that
if $g(a)g(a+d)$ is not equal to this, then it must take the value $0$.
It follows that 
\begin{equation} \label{Z1_find}
\EE(Z_1)\ =\ p^{-2} |S|(|S|-1) (1-\eps)^{-2} W_1\ =\ W_1 + O(p).
\end{equation}
Likewise for the other $Z_i$, we will have that 
$$
\EE(Z_i)\ =\ W_i + O(p).
$$
\bigskip

\noindent {\bf Collecting the two cases together.}
\bigskip

Let $Z_7$ denote the contribution of arithmetic progressions lying entirely in $B$;
that is,
$$
Z_7\ =\ \SUM_{b,b+d,b+2d \in B} f(b)f(b+d)f(b+2d)\ =\ \SUM_{b,b+d,b+2d \in B}
g(b)g(b+d)g(b+2d).
$$
Note that in this case $W_7 = Z_7$.

Putting together our above estimates, and using the fact that 
$$
\LL(g)\ =\ p^{-2} (Z_0 + \cdots + Z_7),
$$
we find that
\begin{eqnarray*}
\EE(\LL(g))\ &=&\ p^{-2}(W_0 + \cdots + W_7 - \eps^2 W_0/2 + O(p) ) \\
&=&\ \LL(f) - \eps^2 p^{-2} W_0/2 + O(1/p).
\end{eqnarray*}
Using Markov's inequality we have
$$
\PP(\LL(g)\ <\ \LL(f) - \eps^2 p^{-2} W_0/4)
\ \geq\ 1\ -\ {\EE(\LL(g)) \over \LL(f) - \eps^2 p^{-2} W_0/4} 
\ >\ \eps^2/8,
$$
since $\LL(f) \geq p^{-2} W_0$.
\bigskip

\noindent {\bf $\EE(g)$ is close to $\EE(f)$ with high probability.}
\bigskip

Before we ``derandomize'' and pass to an instantiation of $g$, we will
need to also show that $\EE(g)$ is close to $\EE(f)$ with high probability.
This can be accomplished in several different ways, though here we will
just use the second moment method:  First, let 
$$
F\ :=\ \SUM_{a \in A} f(a),\ {\rm and\ } G\ :=\ \SUM_{a \in A} g(a).
$$ 
Now, as is easy to show, $F + O(1/p) = \EE(G)$; and so, since 
$\eps\beta > p^{-1/2} \log p$, we have that
\begin{equation} \label{FG}
\PP(|F-G| \geq 2\beta p)\ \leq\ \PP(|G - \EE(G)| \geq \beta p).
\end{equation}
It follows from Chebychev's inequality that this last probability is at most
$$
{ {\rm Var}(G) \over \beta^2 p^2}\ =\ { \EE(G^2) - \EE(G)^2 \over \beta^2 p^2}.
$$
To bound this from above we observe that
$$
\EE(G^2)\ =\ \SUM_{a,b \in A} \EE(g(a)g(b)). 
$$
Now, as a consequence of what we worked out just before 
(\ref{Z1_find}), we have that 
$g(a)$ and $g(b)$ are independent whenever $a \neq b$.  So, 
$$
\EE(G^2)\ =\ \EE(G^2) + O(p),
$$
and it follows that the probability of the right-most event in 
(\ref{FG}) is at most $O(\beta^{-2}/p)$.  It is easy to see that with
probability $1 - O(\beta^{-2}/p)$ we will have
\begin{equation} \label{Egf}
\EE(g)\ \geq\ \EE(f) - 2\beta.
\end{equation}
\bigskip

\noindent {\bf Conclusion of the proof.}
\bigskip

It follows that with probability at least
$$
(1 - O(\beta^{-2}/p))\ +\ \eps^2/8\ -\ 1
$$
we will have that
$$
\EE(g)\ \geq\ \EE(f) - 2\beta\ \ {\rm and\ \ } 
\LL(g)\ \leq\ \LL(f) - \eps^2 p^{-2} W_0/4 + O(1/p).
$$
Using our assumption that
$$
\eps\beta >\ p^{-1/2} \log p,
$$
we have that this probability is positive.  
So, there exists an instantiation of $g$
such that both hold; henceforth, $g$ will no longer be random, but will instead be one
of these instantiations.

By reassigning at most 
$2\beta p$ places $a \in A$ where $g(a) = 0$ to the value $1$, we can
guarantee that $\EE(g) \geq \EE(f)$, and one easily sees that
$$
\LL(g)\ <\ \LL(f) + 2\beta - \eps^2 p^{-2} W_0/4 + O(1/p).
$$
This completes the proof of our proposition.
\hfill $\blacksquare$

\subsection{Proof of Proposition \ref{forcing_prop}} \label{forcing_prop_section}

We have that if we define the new function 
$h_3(n) = h(n)$ at all $n \in \F_p$, 
except for $n=x$ and $n=y$, then 
$$
\LL(h_3)\ =\ \LL(h) + E_1 + \cdots + E_{13},
$$
where if we let $\omega = e^{2 \pi i /p}$, then
\begin{eqnarray*}
E_1\ &=&\ p^{-3} \SUM_a \hat h(a)^2 (h_3(y) - f(y)) \omega^{-2ay}\ =\ p^{-2} (h*h)(2y)
(h_3(y) - h(y)) \\
E_2\ &=&\ p^{-3} \SUM_a \hat h(a)^2 (h_3(x) - h(x)) \omega^{-2ax}\ =\ p^{-2} (h*h)(2x)
(h_3(x) - h(x)) \\
E_3\ &=&\ 2 p^{-3} \SUM_a \hat h(a) \hat h(-2a) (h_3(y) - h(y)) \omega^{ay}\\
&=&\ 2p^{-2} (h*h_2)(-y) (h_3(y) - h(y)) \\
E_4\ &=&\ 2p^{-3} \SUM_a \hat h(a) \hat h(-2a) (h_3(x) - h(x)) \omega^{ax} \\
&=&\ 2p^{-2} (h*h_2)(-x) (h_3(x) - h(x)) \\
E_5\ &=&\ 2p^{-3} \SUM_a \hat h(a) (h_3(y) - h(y))^2 \omega^{-ay}\ =\ 
2p^{-2} (h_3(y) - h(y))^2 h(y) \\
E_6\ &=&\ 2p^{-3} \SUM_a \hat h(a) (h_3(y) - h(y))(h_3(x) - h(x)) \omega^{a(y-2x)} \\
&=&\ 2p^{-2} (h_3(y) - h(y))(h_3(x) - h(x)) h(2x-y) \\
E_7\ &=&\ 2p^{-3} \SUM_a \hat h(a) (h_3(y) - h(y))(h_3(x) - f(x)) \omega^{a(x-2y)} \\
&=&\ 2p^{-2} (h_3(y) - h(y))(h_3(x) - f(x)) h(2y-x) \\
E_8\ &=&\ 2p^{-3} \SUM_a \hat h(a) (h_3(x) - h(x))^2 \omega^{-ax}\ =\ 2p^{-2} (h_3(x) - h(x))^2 
h(x) \\
E_9\ &=&\ p^{-3} \SUM_a \hat h(-2a) (h_3(y) - h(y))^2 \omega^{2ay}\ =\ 
p^{-2}(h_3(y) - h(y))^2 h(y) \\
E_{10}\ &=&\ p^{-3} \SUM_a \hat h(-2a) (h_3(x) - h(x))^2 \omega^{2ax}\ =\ 
p^{-2} (h_3(x) - h(x))^2 h(x) \\
E_{11}\ &=&\ 2p^{-3} \SUM_a \hat h(-2a) (h_3(x) - h(x))(h_3(y) - h(y)) \omega^{a(x+y)} \\
&=&\ 2 p^{-2} (h_3(x) - h(x))(h_3(y) - h(y)) h((x+y)/2) \\
E_{12}\ &=&\ p^{-2} (h_3(y) - h(y))^3 \\
E_{13}\ &=&\ p^{-2} (h_3(x) - h(x))^3.
\end{eqnarray*}
There are actually $6$ more terms that make up the above ``error''; however, all
of these give a contribution of $0$, which is why they were not listed.

So, one sees that
\begin{eqnarray} \label{Lf2}
\LL(h_3)\ &=&\ \LL(h) + p^{-2} (h_3(x) - h(x)) F_h(x) + p^{-2} (h_3(y) - h(y)) 
F_h(y) \nonumber \\
&& \hskip1in + E_5 + \cdots + E_{13}.
\end{eqnarray}
\bigskip

To prove our proposition, all we need to show is that if 
there is a pair $x,y \in \F_p$, $x,y \not \in A$, with 
$$
F_h(x)\ <\ F_h(y),\ {\rm and\ } h(y)\ >\ 0,
$$
then in fact 
$$
h(x)\ =\ 1.
$$
Suppose there were such a pair $x,y$ for which $h(x) < 1$.  
Then, we will show that
$h$ fails to minimize $\LL(h)$ subject to the various constraints:  
Basically, we let 
$$
0\ <\ \eps\ < \min(1-h(x), h(y))
$$
(its exact value will be decided later) and then
we consider the function $h_3$ given by
\begin{eqnarray*}
&& {\rm for\ } n \in \F_p,\ n \neq x,y,\ {\rm we\ set\ } h_3(n)\ =\ h(n);\ 
{\rm and,}\\ 
&& h_3(x)\ =\ h(x) + \eps,\ h_3(y)\ =\ h(y) - \eps.
\end{eqnarray*}

From our formula (\ref{Lf2}), we easily deduce that
$$
\LL(h_3)\ \leq\ \LL(h) + \eps p^{-2} (F_h(x) - F_h(y)) - 
O(\eps^2 p^{-2}).
$$
Clearly, if we take $\eps > 0$ small enough, we will get 
$$
\LL(h_3)\ <\ \LL(h),
$$
which contradicts the minimality of $h$.  We conclude, therefore, that 
$h(x) = 1$, as claimed.
\hfill $\blacksquare$

\end{document}